# Harvesting Fisheries Management Strategies With Modified Effort Function

**Lev V. Idels***
Department of Mathematics, Malaspina University- College,
Bld. 360/304,  900 Fifth Street,  Nanaimo, BC,  V9R 5S5
Tel: (250) 753-3245 ext. 2429
Fax: (250) 740-6482
email: lidels@shaw.ca
*Corresponding author

**Mei Wang**
Institute of Applied Mathematics,
University of British Columbia, BC, Canada
E-mail: maggie@math.ubc.ca

**Biographical notes:**

**Dr. Lev Idels** is a University-College Professor at the Department of Mathematics at Malaspina University-College in Nanaimo, BC. His research interests have revolved around the qualitative analysis of nonlinear differential equations and nonlinear differential equations with delay and their applications in biosciences  He is teaching a great variety of Math Courses, including Mathematical Modeling. In 2005 Dr. Lev was awarded an NSERC Discovery Grant for the research in Mathematical Biology.

**M. Wang** is a PhD student in Institute of Applied Mathematics of University of British Columbia, Canada. Her current research interests include numerical solution of Partial Differential Equations (PDE), combustion simulation and all kinds of science applications of PDE and ODE (ordinary differential equations).

**Abstract**

In traditional harvesting model, a fishing effort, *E,* is defined by the fishing intensity and does not address the inverse effect of fish abundance on the fishing effort. In this paper, based on a canonical differential equation model, we developed a new fishing effort model which relies on the density effect of fish population. We obtained new differential equations to describe certain standard Fisheries management strategies. This study concludes that a control parameter $\beta$ *(* the magnitude of the effect of the fish population size on the fishing effort function *E),* changes not only the rate at which the population goes to equilibrium, but also the equilibrium values.

To examine systematically the consequences of different harvesting strategies, we used numerical simulations and qualitative analysis of six fishery strategies, e.g., proportional harvesting, threshold harvesting, proportional threshold harvesting, and seasonal and rotational harvestinng.

## 1.   Introduction

A major current focus in fishery management is how best to ensure harvesting sustainability [1-3,5,7,13]. Clearly the object of the management is to devise harvesting strategies that will not drive species to extinction. Therefore, the notion of persistence and extinction times of the populations, as well as a precautionary harvesting policy, is always critical. A control variable of every fishery management is the fishing effort [2,8], which is defined as a measure of the intensity of fishing operations.

The Schaefer fishing model [4,5,7,14] takes the form:

$$\frac{dN}{dt} = rN(t)\left[1 - \frac{N(t)}{K}\right] - Y(t) \tag{1}$$

where *N* is the population biomass of fish at time *t*, *r* is the intrinsic rate of growth of the population, *K* is the carrying capacity, and we assume that *r*≥0 and *K*>0 are constants. The harvest function is defined as

$$Y(t) = qN(t)E \tag{2}$$

Here *q*≥0 is the catchability coefficient, defined as the fraction of the population fished by a unit of effort. *E*≥0 is the fishing effort, the intensity of the human activities to extract the fish. In general, fishing effort is regulated by quotas, trip limits and gear restrictions.

Equation (2) implies that harvest per unit effort is a function of the size of the population.

$$\frac{Y(t)}{E} = qN(t) \tag{3}$$

If the price of fish responds to the quantity of the harvest, a greater harvest would induce a lower price of harvest, and vice versa. If we assume that the market price of the harvest motivates changes in fishing effort, a lower price (or a larger population) induces less fishing effort, and vice versa. In traditional fishery models [2,5,8,9], fishing effort $E$ is simply expressed as a function of time $E=E(t)$, which does not address the inverse effect of fish abundance on the fishing effort [11] (higher density of fish, less effort to catch unit harvest). Under the above assumptions, it is more appropriate to express $E$ as a function of population dynamics. For example,

$$E(t,N) = \alpha(t) - \beta(t)\frac{1}{N}\frac{dN}{dt} \quad (4)$$

where $\alpha \geq 0$ and $\beta \geq 0$ are continuous functions of $t$. If we substitute (4) back into Eq. (2) we have

$$Y(t) = qN(\alpha(t) - \beta(t)\frac{1}{N}\frac{dN}{dt}) \quad (5)$$

We substitute (5) into Eq. (1):

$$\frac{dN}{dt} = rN[1-N/K] - qN(\alpha(t) - \beta(t)\frac{1}{N}\frac{dN}{dt}). \quad (6)$$

An often debated question [3, 4,11] is whether and when seasonal harvesting strategies are effective in fishery management. Traditional management strategies for continuous harvest may lead to serious mismanagement (e.g., extinction of fish population). There has been growing interest in rotational use of fishing grounds [9,10,13]. There is some evidence that rotational use of fishing grounds slightly increases both yield- and biomass-per-recruit [13].

For our numerical simulations three practical data sets were chosen from [6,11]:

data1    $r$=0.3; $\alpha q$=0.24/0.36/0.42/1.00 and $E$=1.15 for the max sustainable yield;
data2    $r$=0.3; $\alpha q$=0.57/1.58/1.70/1.81 and $E$=0.8 for the max sustainable yield;
data3    $r$=0.3; $\alpha q$=0.04/0.18/0.49/0.61/0.72/1.0 and $E$=0.25 for the max sustainable yield

## 2. Standard fishery strategies models, qualitative analysis and numerical solutions

### 2.1 Constant harvesting

Constant harvesting removes a fixed number of fish each year. Furthermore we assume fishermen have perfect information of where the fish are

$$\frac{dN}{dt} = rN[1-N/K] - qE \quad (7)$$

### 2.2 Proportional harvesting

A constant fraction of fish is removed each year:

$$\frac{dN}{dt} = rN[1-N/K] - \lambda qN(\alpha(t) - \beta(t)\frac{1}{N}\frac{dN}{dt}) \quad (8)$$

Where $\lambda$ is the proportional rate and $0 \leq \lambda \leq 1$ is a constant. This equation can be simplified to:

$$\frac{dN}{dt} = \frac{r}{1-\lambda q\beta} N[1-N/K] - \frac{\lambda q\alpha}{1-\lambda q\beta} N \qquad (9)$$

where $\lambda q\beta \neq 1$.

- Qualitative analysis

Equation (9) has an explicit solution

$$N(t) = \frac{(1-\lambda q\alpha)K}{1+C\exp\left(\frac{-r(1-\lambda q\alpha)}{1-\lambda q\beta} t\right)}$$

The two equilibrium point for this system is $N_u^* = 0$ and $N_s^* = (1-\frac{\lambda q\alpha}{r})K$, where $N_u^*$ and $N_s^*$ represent the unstable and stable equilibrium point separately. In order to have a stable point with positive fish population $N_s^* > 0$ or

$$\lambda q\alpha < r \qquad (10)$$

Example: Let $r=0.5$, $\lambda=0.5$, $q=0.8$, $\beta=1.0$, $\alpha=1.0$, $K=1$, and $\lambda q\alpha=0.4<0.5$. Figure 1 is the phase portrait of this system. The equilibrium point is the intersection of $\frac{r}{1-\lambda q\beta} N(1-N/K)$ and $\frac{\lambda q\alpha}{1-\lambda q\beta} N$ which is $N_s^* = 0.2$ in this case. As the slope of $\frac{\lambda q\alpha}{1-\lambda q\beta} N$ grows, the value of positive $N_s^*$ is decreased and the fish population becomes extinct. The system undergoes a transcritical bifurcation as $\lambda q\alpha$ increases.

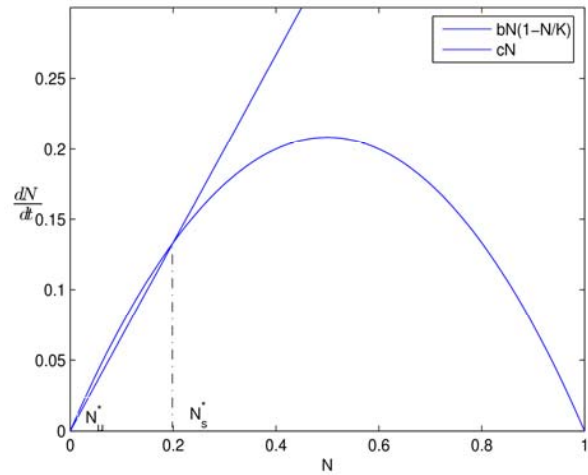

Figure 1: Phase portrait of proportional harvesting. $r=0.5$, $\lambda=0.5$, $q\alpha=0.8$, $\beta=1.0$ and $K=1$.

If $\lambda=0.5$, that satisfies the condition of Eq.(10), all solutions with different initial conditions will approach $N_s^*=0.2$. As $\lambda$ increases (or an increase in $\alpha q$), the population is overfished and goes to extinction eventually.

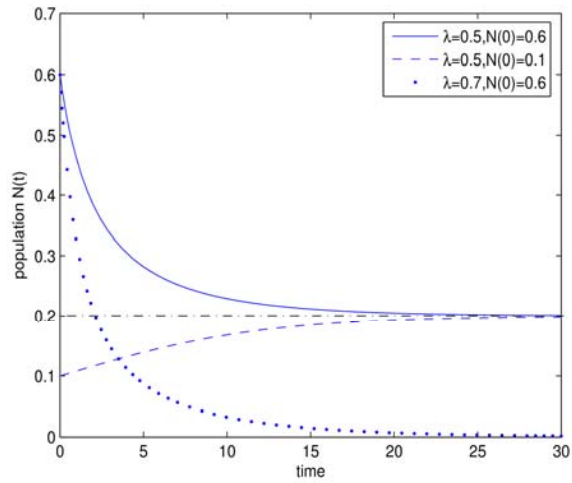

Figure 1: Population biomass dynamics with proportional harvesting. r=0.5, α q=0.8, β=1, K=1

- Sustainable yield and dynamics of $\beta$

The sustainable yield at equilibrium is:

$$Y = \lambda q N_s^* \alpha = \lambda q \alpha K (1 - \frac{\lambda q \alpha}{r}) \qquad (11)$$

Let $\dfrac{dY}{d(\lambda q \alpha)} = 0$, then $\lambda q \alpha = \dfrac{r}{2}$ And

$$Y_{max} = \frac{rK}{4} \qquad (12)$$

In Figure 3, as the proportion rate increase, the value of Y will increase but after a certain limit, the fish population will become extinct.

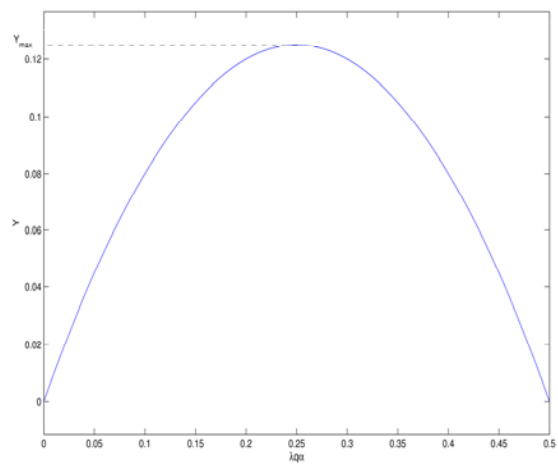

Figure 3: Sustainable yield of proportional harvesting. r=0.5.

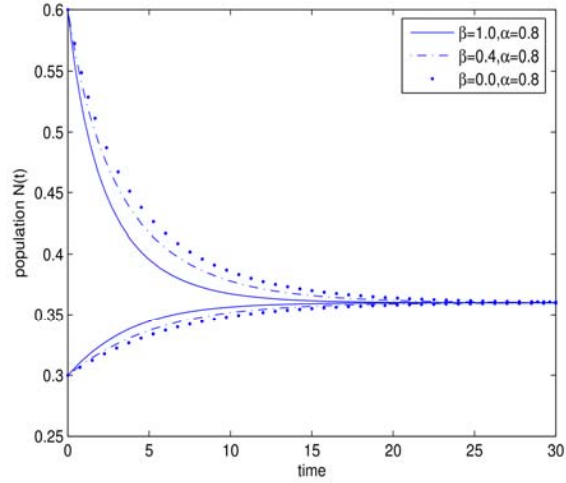

Figure 4: Dynamics of β, r=0.5, β=0.5, q=0.8.

Figure 4 shows a time series of $N(t)$ with three $\beta$ values which range from above or below $\alpha$. As $\beta$ increases, $N(t)$ converges to the equilibrium solution faster for any value of $N(0)$. That is, as the value of β increases, $E$ will be affected more by the change of fish population, thus it will cause the $N(t)$ to change and converge faster.

## 2.3  Restricted proportional harvesting

An upper limit of harvest is introduced on proportional harvesting:

$$\frac{dN}{dt} = \begin{cases} rN[1-N/K] - \lambda qN(\alpha(t) - \beta(t)\frac{1}{N}\frac{dN}{dt}) & \text{if } Y(t) \leq Y_{\text{limit}} \\ rN[1-N/K] - Y_{\text{limit}} & \text{if } Y(t) > Y_{\text{limit}} \end{cases} \quad (13)$$

where $Y_{\text{limit}} \geq 0$. As long as $0 \leq Y_{\text{limit}} \leq \frac{rK}{4}$, when $Y(t) \leq Y_{\text{limit}}$, the population will behave the same as strategy 2. As $Y(t) > Y_{\text{limit}}$ we will have two positive equilibrium points shown in Figure 5. One thing to notice is that if the initial value $N(0) < N_u^*$ then the fish population will become extinct. Therefore, $Y_{\text{limit}}$ should be constrained by N(0), i.e., $Y_{\text{limit}} < rN(0)(1 - \frac{N(0)}{K})$.

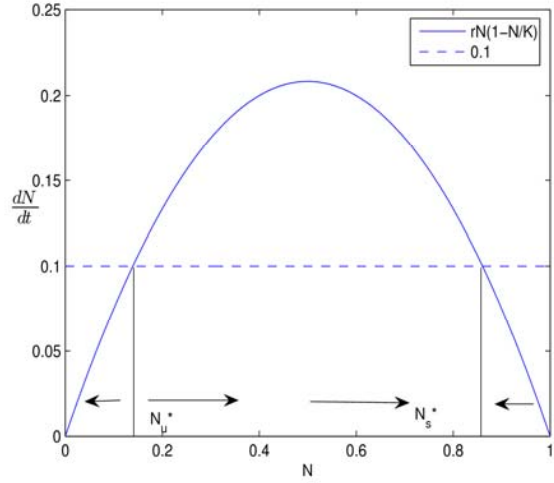

Figure 5: Phase portrait of restricted proportional harvesting.

## 2.4 Proportional threshold harvesting

In proportional threshold harvesting, only a fixed proportion of the fish above the threshold is harvested.

$$\frac{dN}{dt} = \begin{cases} rN[1-N/K] - \lambda(N-N_{thre})qN(\alpha(t) - \beta(t)\frac{1}{N}\frac{dN}{dt}) & \text{if } N(t) > N_{thre} \\ rN[1-N/K] & \text{if } N(t) \leq N_{thre} \end{cases} \quad (14)$$

where $N_{thre} \geq 0$. In this case the excess stock is harvested. This policy provides the stock with a measure of protection from overexploitation. Eq. (14) can be simplified as:

$$\frac{dN}{dt} = \frac{rN(1-N/K) - \lambda q\alpha(N-N_{thre})}{1 - \lambda q\beta(1-\frac{N_{thre}}{N})} \quad (15)$$

$$\text{where } N_{thre} \neq N(1-\frac{1}{\lambda\beta q})$$

• Qualitative analysis

From the analytical solution of this equation for $N_{thre} > 0$ we have $N_u^* \leq 0$ and $N_s^* \geq 0$. As the slope of the straight line or $N_{thre}$ decreases, $N_s^*$ increases and had been always less than $K$. Because $N(0) \geq 0$, all the solutions will converge to $N_s^* \geq 0$. Comparing with strategy 2 with the same set of parameters, $N_s^*$ increases but $Y$ decreases in this strategy.

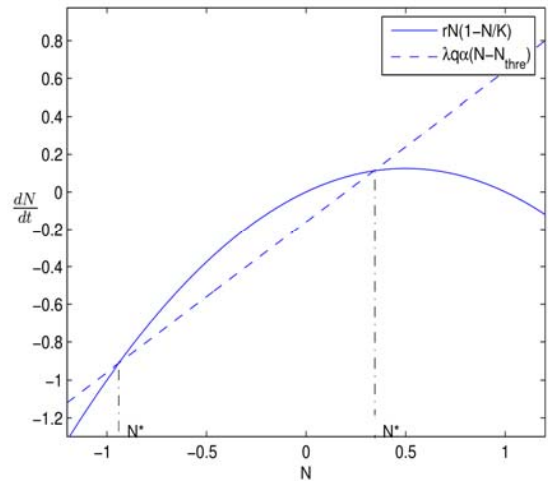

Figure 6: Phase portrait of proportional threshold harvesting. r=0.5, λ=0.5, α q=0.8, β=1.0, K=1, Nthre=0.2

- Sustainable yield and effects of β

The sustainable yield at equilibrium is:

$$Y = \lambda q N_s^* \alpha (N_s^* - N_{thre})$$

and the total yield will increase proportionally to λqα.

Figure 7 shows the change of $N(t)$ for different values of β which vary from above or below α. It is interesting to note that as β increases, $N(t)$ will converge to the $N_s^*$ faster if $N(0)$ is above the equilibrium solution. However, if $N(0)$ is below $N_s^*$, an increase in β will cause $N(t)$ to converge slightly slower.

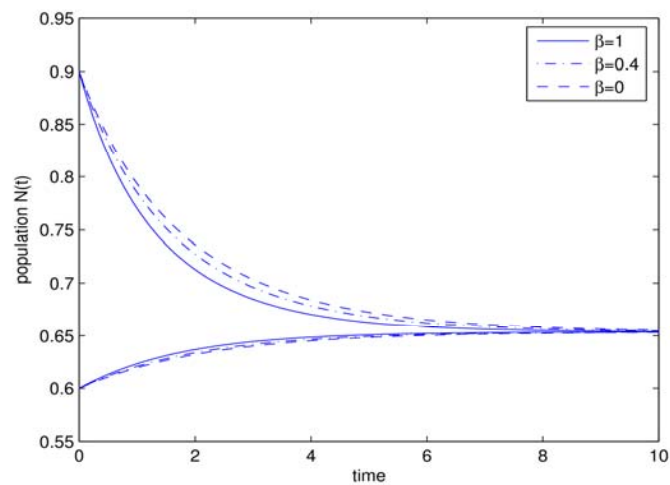

Figure 7. Effects of β in proportional threshold harvesting $r = 0.1$; $\lambda = 0.5$; q=0.4; K = 0.5; $N_{thre} = 0.3$.

## 2.5 Seasonal (periodic) harvesting

$$\frac{dN}{dt} = rN[1-N/K] - \lambda(t)qN(\alpha(t) - \beta(t)\frac{1}{N}\frac{dN}{dt}) \quad (17)$$

where $\lambda(t)$ is a periodic function of time with the period of 1 year. This system has the same format as case 2, thus the limitation is $\lambda q\alpha < r$. However, because $\lambda$ is a periodic function and varies from season to season, the fish won't become extinct during fishing time and, if in some season we stop fishing, the amount of fish might be able to increase again.

1. First example:

$$\lambda(t) = \begin{cases} 0.5\sin(\pi\frac{(t-n-t_{start})}{H}) & f \ n+t_{start} < t < n+t_{start} + H, \ n=0,1,2... \\ 0 & other \end{cases}$$

where $t_{start}$ is the harvest starting time within one year. The first test is based on $t_{start}=0.25$ and $H=0.25$ which means harvest in summer season only.

In a long run, the population is shown in Figure 8:

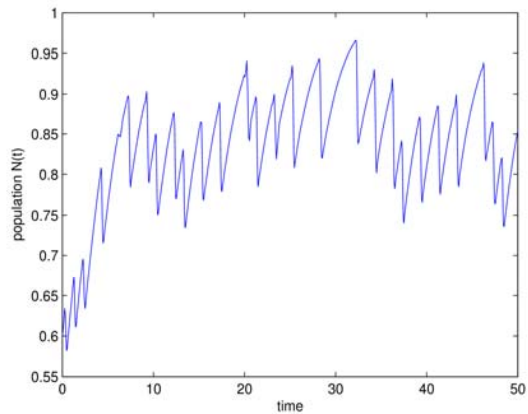

Figure 8: Long-time behavior of population in seasonal harvesting $r=1, \lambda=0.5, q\ \alpha=0.4$.

A solution of this system depends on the magnitude of $\lambda$ and how it varies during different seasons.

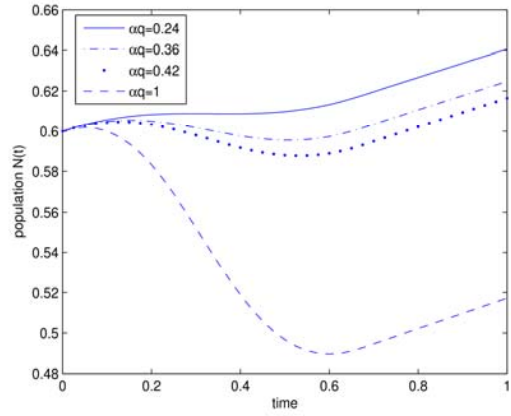

Figure 9: Seasonal harvesting one year solution with data1

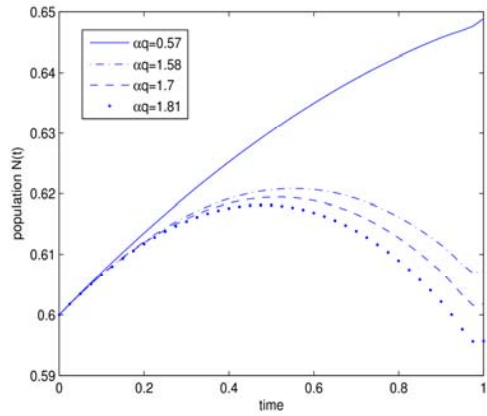

Figure 10: Seasonal harvesting one year solution with data2

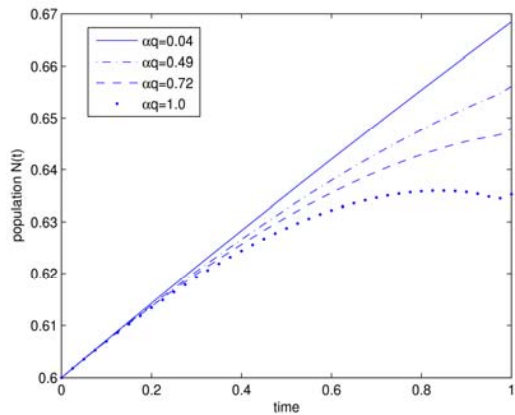

Figure11. Seasonal harvesting one year solution with data3

On Figure 12 we repeat Figure 11 to show a long-term behaviour of the system

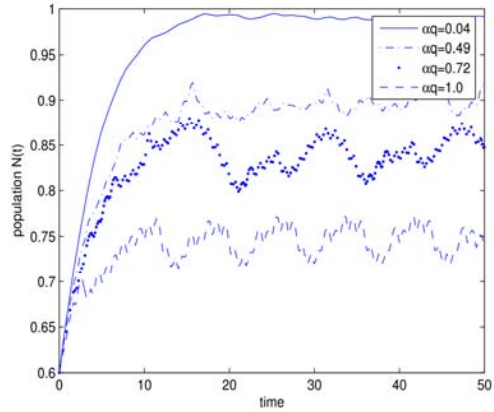

Figure 12: Seasonal harvesting with data3

2. Second example:
$$\lambda(t) = \begin{cases} 0.5 & \text{if } Hn + bn < t < H(n+1) + bn;\ n = 0, 1, 2... \\ 0 & H(n+1) + bn < t < H(n+1) + b(n+1);\ n = 0, 1, 2... \end{cases}$$

We applied this example to data1 in Figure 13.

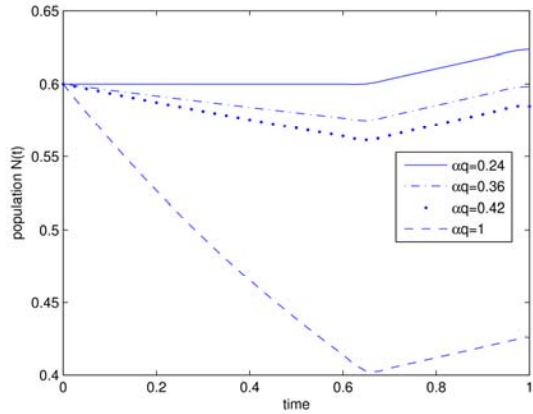

Figure 13: One year solution with data1

3. Third example:
$$\lambda(t) = \frac{1 + \sin(2\pi t)}{4} \tag{18}$$

whose maximum is $\lambda(0.25+n)=0.5$ and minimum $\lambda(0.75+n)=0, n=0,1,2....$

Figure 1s 13-16 are the numerical solutions with data1, data2 and data3.

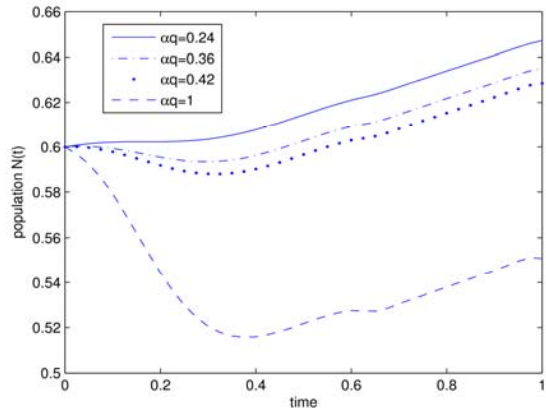

Figure 14: One year solution with data1

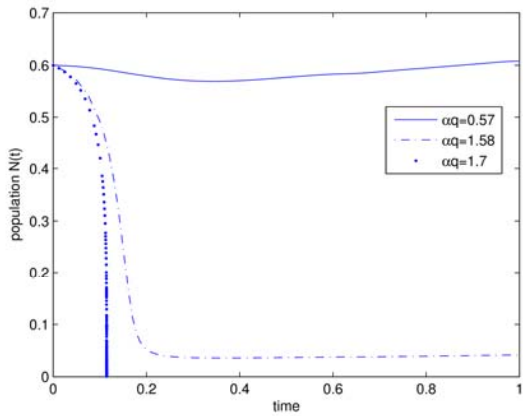

Figure 15: One year solution with data2

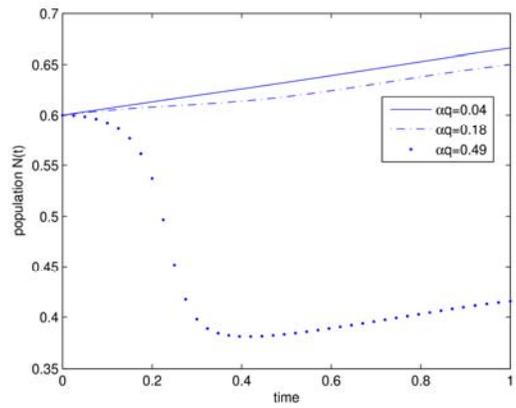

Figure 16: One year solution with data3

4. The effects of β

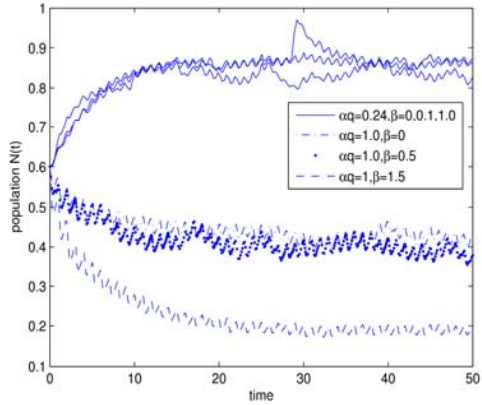

Figure 17: Dynamics of seasonal harvesting with data1.

Figure 17 shows that β changes the $N_s^*$ only if $N_s^*$ is below $N(0)$, specifically, as β decreases, $N_s^*$ increases.

In another test case shown in Figure 18, as β changes the rate of $N(t)$ is approaching $N_s^*$. However, it does not change the value of $N_s^*$ itself.

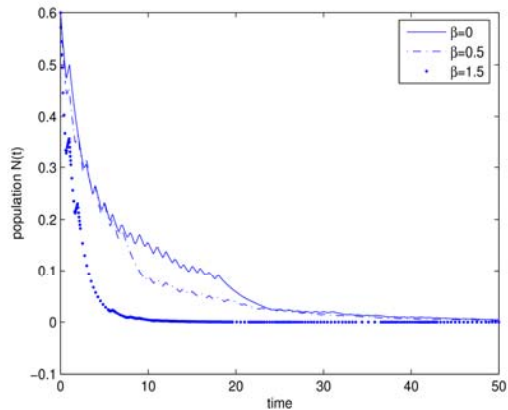

Figure 18: Short time behavior and the dynamics of β in seasonal harvesting with data3.

Comparing Figure 17 and Figure 18, we can see that the effect of β depends on the particular character of $\lambda(t)$ itself.

Comparing the seasonal harvesting with strategy 2 and 4 which was tested with the same set of data, we can see that strategy 6 behaves similarly to strategy 2 depending on $\lambda(t)$. Also due to the periodic nature of λ, it tends to have a higher $N_s^*$ than strategy 2. Although strategy 4 has no worry of extinction because of the existence of $N_{thre}$, it has the cost of decrease of $Y$.

## 2.6 Rotational harvesting

In practice, this strategy is usually executed by dividing the field into different areas and closes some of them rotationally [7,11]. Assume there is no immigration and migration of fish, and consider just total area of rotation as one dynamic system, then the system will behave almost the same way as periodic harvesting. The numerical analysis will be almost the same as periodic harvesting, the difference is just that the closure period will likely be longer and the fish population has more time to recover from the harvest season. Therefore it is possible to increase the maximum yield within a certain period of time. Therefore, only one numerical simulation example is listed here. Consider the third example in periodic harvesting with data1 and harvesting for one year and then close for two years.

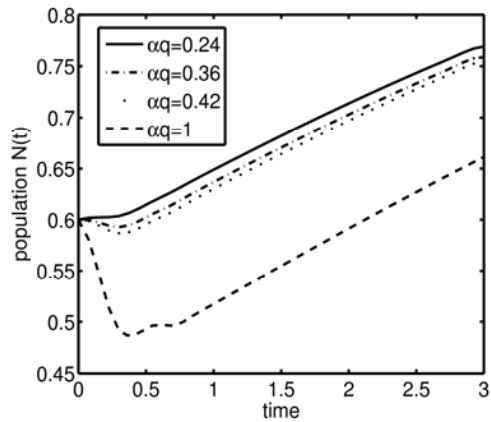

Figure19: Short time behavior in rotational harvesting

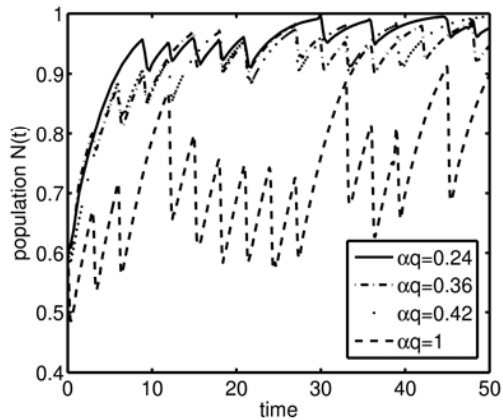

Figure20: Long-time behavior in rotational harvesting

# 3. Conclusions

We developed a new fishing effort model which relies on the density effect of fish. We study the consequences of harvesting with five fishery strategies. This study concludes a control parameter $\beta$, which defines the magnitude of the effect of the fish population size on *E,* changes not only the rate at which the population goes to equilibrium, but also the equilibrium values.

It is shown that β plays a significant role for the seasonal harvesting. Comparing the seasonal harvesting with strategy 2 and 4, that was tested with the same set of data, we can see that strategy 5 has a similar behavior as strategy 2 depending on λ(*t*). Also due to the period behavior of λ, it tends to have a higher $N_s^*$ than strategy 2. Although strategy 5 prevents fish population from extinction because of the threshold population size $N_{thre}$.

It follows from our analysis that the rotational use of fishing grounds increases both yield- and biomass-per-recruit, while still keeps the fish population sustainable.

While the rotational harvesting behavior is almost the same as the case of periodic harvesting, the difference appears when q=1: while the continuous periodic harvesting will cause the fish to extinguish, it is not the case for rotational harvesting, which clearly fluctuate a lot with a period of 3 but managed to remain sustainable and fishes are not extinguished. It proves that the rotational harvest strategy is good in a sense that it can yield a larger amount of fish in a certain period of time.

Table 1 summarizes the results for different strategies.

| Strategy | Limitations | Dynamics |
|---|---|---|
| Proportional harvesting (Strategy 2) | λ*q*α<*r* | As the slope of $\frac{\lambda q \alpha}{1-\lambda q \beta}N$ grows, the value of positive $N_s^*$ is decreased and severe over-fishing happens and the fish approach extinction. The system undergoes a transcritical bifurcation as λ*q*α increases. Increase in β will stabilize the population faster |
| Proportional threshold harvesting (Strategy 4) | *N*(0)>0 | Increase in β will converge *N*(*t*) to $N_s^*$ faster if *N*(0) is above the equilibrium solution. However, if *N*(0) is below $N_s^*$, increase in β will cause *N*(*t*) to converge to $N_s^*$ slower. |
| Seasonal harvesting (Strategy 5) | Dependence on λ(*t*) | Also due to the periodic nature of λ, it tends to have a higher $N_s^*$ than strategy 2. Although strategy 5 has no worry of extinction because of the existence of $N_{thre}$, it has the cost of a decrease in *Y*. Changes the values of $N_s^*$ |

Table 1: Different strategies comparison